\documentclass[a4paper,12pt,leqno]{amsart}

\usepackage{amssymb,hyperref}%,showkeys}
\usepackage[latin1]{inputenc}
\swapnumbers
\theoremstyle{plain}

\newtheorem{theo}{Theorem}[section]
\newtheorem{prop}[theo]{Proposition}
\newtheorem{cor}[theo]{Corollary}
\newtheorem{lem}[theo]{Lemma}

\theoremstyle{definition}
\newtheorem{defi}[theo]{Definition}
\newtheorem{notn}[theo]{Notation}

\newtheorem{question}[theo]{Question}
\theoremstyle{remark}
\newtheorem{rem}[theo]{Remark}
\newtheorem{ex}[theo]{Example}
\numberwithin{equation}{section}

\newcommand{\thismonth}{\ifcase\month\or
  January\or February\or March\or April\or May\or June\or July\or
  August\or September\or October\or November\or December\fi
  \space\number\year}
\DeclareMathAlphabet{\mathrmsl}{OT1}{cmr}{m}{sl}

\newcommand{\oper}[3][n]{\newcommand{#2}{\mathop{\mathrm{#3}}%
\ifx n#1\nolimits\else\limits\fi} }
\newcommand{\rsoper}[3][n]{\newcommand{#2}{\mathop{\mathrmsl{#3}}%
\ifx n#1\nolimits\else\limits\fi} }
\newcommand{\TM}{{\mathbb T}}
\newcommand{\RM}{{\mathbb R}}
\newcommand{\CM}{{\mathbb C}}
\newcommand{\ra}{\rangle}
\newcommand{\la}{\langle}
\renewcommand{\div}{\operatorname{div}}

\newcommand{\dbarb}{\bar{\partial}_{b}}

\newcommand{\Scal}{\operatorname{Scal}}
\newcommand{\Ric}{\operatorname{Ric}}

\newcommand{\sym}{\operatorname{Sym}_2^*H^{(1,0)}}

\newcommand{\cerc}{{\mathbb S}}

\renewcommand{\exp}{\operatorname{e}}

\newcommand{\cnabla}{{\mathfrak D}}
\renewcommand{\Re}{\operatorname{Re}}
\renewcommand{\Im}{\operatorname{Im}}

\newcommand{\kn}{\,\protect{\bigcirc\hspace{-10pt}\wedge}\,}
\newcommand{\proofof}[1]{\end{#1}\begin{proof}}

\def\endproof{\hfill\qed}
\newcounter{mnotecount}[section]
\renewcommand{\themnotecount}{\thesection.\arabic{mnotecount}}
\newcommand{\mnote}[1]%{}
{\protect{\stepcounter{mnotecount}}$^{\mbox{\footnotesize  $%\!\!\!\!\!\!\,
      \bullet$\themnotecount}}$ \marginpar{\raggedright\tiny\em
    $\!\!\!\!\!\!\,\bullet$\themnotecount: #1} }

\begin{document}

\title[Cartan connection and CR geometry]{The canonical Cartan bundle and connection\\ in CR geometry}
\author{Marc Herzlich}
\address{Institut de Math\'ematiques et Mod\'elisation de Montpellier\\ 
UMR 5149 CNRS -- Universit\'e Mont\-pellier~II\\ France}
\email{herzlich@math.univ-montp2.fr}

\begin{abstract}
We give a differential geometric description of the canonical Cartan (or tractor) bundle and connection
in CR geometry, thus offering a direct, alternative, definition to the usual abstract approach.
\end{abstract}

\maketitle

\bigskip
\section{Introduction.}

The goal of this paper is to present a simple differential geometric description
of the Cartan connection naturally associated to a strictly pseudoconvex (integrable)
CR manifold.

Cartan connections have proved to be useful tools in CR geometry, and more generally
in the so-called \emph{parabolic geometries}, for exhibiting invariants or invariant 
differential operators, see \emph{e.g.} \cite{gover-graham}. Given a semisimple Lie group $G$ 
and a parabolic subgroup $P$, the Cartan connection for the parabolic 
geometry modelled on the homogeneous space $G/P$ (CR geometry being an example thereof) is 
usually described as a $P$-equivariant $1$-form $\omega$ on a principal $P$-bundle $\mathcal{G}$ 
with values in the Lie algebra $\mathfrak{g}$ that establishes 
an isomorphism between $T_s\mathcal{G}$ and $\mathfrak{g}$ at each point $s$ in $\mathcal{G}$.
Its construction usually follows from elaborate Lie algebra considerations \cite{cap-schichl}.

Given some faithful representation of $G$, one can construct a vector bundle endowed with a linear
covariant derivative from $({\mathcal G},\omega)$. In CR geometry, one such choice is the standard
representaton of $G=SU(1,n+1)$ and the resulting vector bundle $\mathbb{T}$ is usually called 
the \emph{Cartan bundle} or \emph{tractor bundle}. Exhibiting the bundle $\mathbb{T}$ and its
covariant derivative provides another (slightly simpler in its formulation) definition 
of the structure of CR geometry as a Cartan geometry. However this is slightly deceptive as 
the known constructions of the bundle $\mathbb{T}$ either involve the whole parabolic geometry 
machinery as above or is done by first giving some \emph{ad hoc} definition of $\mathbb{T}$ 
when some choice of a contact form underlying the contact structure has been made, and then 
checking that the given objects are in fact independent of the contact form. 

The goal of this note is to give an intrinsic definition of $\mathbb{T}$ and its canonical 
connection in terms of a down-to-earth viewpoint on CR geometry, somehow in the spirit of previous
works on conformal geometry due to Gauduchon \cite{pg-cartan} and Bailey, Eastwood, and Gover \cite{beg}. 
More precisely, we start from 
the simplest definition for CR geometry (a codimension $1$ contact distribution in the tangent 
bundle with a compatible complex structure) and we build $\mathbb{T}$ directly from this data, 
without any reference to the general setting of parabolic geometries. Part of the constructions 
here are implicit in many works in CR geometry, especially those of \v{C}ap 
\cite{cap-survey,cap-twistor}, \v{C}ap and Slov\'ak \cite{cap-slovak}, 
and Gover and Graham \cite{gover-graham}. We
also mention that of Farris \cite{farris}, for its use of closed sections of the canonical CR bundle
as a building block for the theory. The 
author however hopes that the explicit definitions detailed below might be useful for a better 
understanding of CR geometry and its canonical Cartan connection.
 
\bigskip

\section{CR geometry and the Cartan vector bundle}

We recall that a strictly pseudoconvex CR manifold is a manifold $M$ (always taken here orientable)
of odd dimension $2n+1$, endowed with a contact structure (a maximally non-integrable distribution of 
hyperplanes $H$ in the tangent bundle) together with a compatible (formally integrable) 
complex structure $J$ on $H$, such
that for any choice of contact form $\eta$, 
$\gamma=d\eta(\cdot, J\cdot)$ is a positive definite quadratic form
at each point on $H$. The existence of $J$ implies the existence of natural sub-bundles of forms of
type $(p,q)$ on the complexified bundle of forms of every degree (note however that, contrarily to the 
usual complex case, the bundles of $(1,0)$- and $(0,1)$-forms do have a non-trivial intersection: the 
$1$-forms whose kernel is precisely $H$ at each point).

\smallskip

Our main object of interest will be the \emph{canonical CR line bundle} $K$, which is the bundle of 
complex $(n+1)$-forms on $M$ of maximal holomorphic type $(n,0)$. Elements of $K$ are $(n+1)$-forms
$\varphi$ such that for any elements $h_1, \dots, h_n$ in $H$,
$$ \varphi(Jh_1,\dots,Jh_n,\cdot) = i^n\, \varphi(h_1,\dots,h_n,\cdot) .$$
Given a contact form $\eta$, one can define a unique \emph{Reeb vector field} $\xi$ transverse to
$H$ by
$$ \eta(\xi) = 1 , \ \ d\eta(\xi,\cdot) =0 .$$
Given a local complex basis $Z_1,\dots,Z_n$ of the space $H^{(1,0)}$ of vectors in $H\otimes\CM$ of 
type $(1,0)$, and the dual basis $\theta^1,\dots,\theta^n$ of horizontal forms of type $(1,0)$, 
the bundle $K$ is generated at each point by $\eta\wedge\theta^1\wedge\dots\wedge\theta^n$.

The canonical CR line bundle has a natural structure of a CR holomorphic line bundle, in the sense
that it has a natural $\dbarb$-operator. More precisely, for any section $\varphi$ of $K$ and any 
vector $Z$ of type $(1,0)$,
$$ \dbarb\varphi (\bar{Z}, \cdots ) =  d\varphi (\bar{Z}, \cdots ) . $$
In case the CR structure of $M$ is induced from an embedding in $\CM^{2n+2}$, the CR holomorphic bundle 
$K$ is nothing but the restriction of the canonical bundle of the complex vector space to $M$ with its 
usual structure of holomorphic line bundle. 

\begin{defi}
The line bundle of $1$-densities $L$ is a $(n+2)$-th root of the dual of the 
canonical CR line bundle $K^*$. 
It has a natural CR holomorphic structure induced from that of $K$.
\end{defi}

In case such a root does not exist, it can at least be defined locally, so that all our 
subsequent constructions keep their meaning in the neighbourhood of a point. This is for example 
enough to study local CR geometry (\emph{e.g.} seeking CR invariant quantities or differential 
operators...)

A possible justification for the name ``bundle of $1$-densities'' is given by the following Lemma, 
which will be of constant use in our work.

\begin{lem} Let $\ell_x$ be a non-zero element of $L$ at a point $x$ of $M$. Then, there exists a unique 
oriented contact form $\theta_{\ell_x}$ at $x$ 
such that $\ell_x$ is of unit norm with respect to the induced pseudo-hermitan structure.
\end{lem}

\begin{proof}
Take an arbitrary coframe $\theta,\theta^1,\ldots,\theta^n$ as above and let $\ell_0$ an element
of $L$ such that $(\ell_0)^{\otimes [-(n+2)]} = \theta\wedge\theta^1\wedge\ldots\wedge\theta^n$.
Write $\ell_x=\exp^{\rho+i\varphi}\ell_0$, then we seek a contact form $\exp^{-2u}\theta$ such that
$$ 1 = |\ell_x|_{\exp^{-2u}\theta} = \exp^{\rho} |\ell_0|_{\exp^{-2u}\theta} =  \exp^{\rho - u} .$$
This leads to the unique choice $u = \rho$. 
\end{proof}

\smallskip

When a local section $\ell$ of $L$ is given, the associated contact form will be called the 
\emph{normalized contact form}. Equivalently, $\ell$ is said to be \emph{volume normalized} 
\cite{farris,lee-fefferman}.

\smallskip

\begin{defi}
The \emph{Cartan vector bundle} $\mathbb{T}$ is the bundle of $1$-jets of CR-holomor\-phic 
sections of $L$.
\end{defi}
 
\smallskip

CR-holomorphic sections of $K$ play an important role in CR geometry. For instance, if 
$M$ is embedded in $\CM^{n+1}$, the natural volume form $dz^1\wedge\dots\wedge dz^{n+1}$ induces
such a section of $K$. F. Farris has used in \cite{farris} the $2$-jet geometry of the 
$\bar\partial_b$-equation for sections of $K$ to give an explicit construction of the Fefferman metric
on the circle bundle $K/\RM_+^*$, see also \cite{lee-fefferman}. As we will see, the jet geometry
of this equation not only yields a nice description of the Fefferman metric but of the full Cartan
bundle.

\smallskip

The Cartan vector bundle is then an $(n+2)$-dimensional complex vector bundle. The choice of a 
connection $\nabla$ on $L$ enables to split the bundle of $1$-jets of sections of $L$ as $J^1L = T^*M 
\otimes L \,\oplus\, L$. If the connection is of type $(0,1)$, \emph{i.e} $\nabla_{|H}^{0,1}
= \bar\partial_b$,
this induces a decomposition of $\TM$, and if moreover some vector field 
$\xi$, transverse to $H$, is given, this may be refined into
\begin{equation}\label{eq:Tfirstdec}
\TM \ = \ L \,\oplus\, (H^*)^{(1,0)} \otimes L \,\oplus\, L  \ \ \ \ ; 
\end{equation}
For instance, a specific choice of contact form $\eta$ gives rise to the Reeb vector field $\xi$ 
defined above and to the so-called Tanaka-Webster connection $\nabla$ \cite{tanaka,webster}, thus 
to the decomposition \eqref{eq:Tfirstdec} above ; for an holonomic element of $\TM$, \emph{i.e.} a 
section of the bundle of $1$-jets of $L$ that is indeed the $1$-jet of some holomorphic section 
$\ell$ of $L$, this simply amounts to identifying the $1$-jet of $\ell$ at a point to the triple
$\left( \ell , \nabla^{1,0}\ell , \nabla_{\xi}\ell  \right)$.  

\bigskip

\section{Weyl structures for CR geometry}

In this section, we shall look a bit further at decompositions of $\TM$ induced by a choice of
a linear connection on $L$. Although not strictly necessary for our purposes, its goal is twofold:
to make contact with explicit expressions that have been given in previous works for the
Cartan vector bundle \cite{cap-slovak,gover-graham}; and to write out a few formulas that will
be useful in our future treatment of the Cartan connection.

\smallskip

We have made use above of the Tanaka-Webster connection for trivialising the bundle $\TM$. A quick 
look at the transformation formulas for the Tanaka-Webster connection when conformally rescaling 
a contact form $\eta$ to $\exp^{2f}\eta$ shows that the last component in \eqref{eq:Tfirstdec}
 (\emph{i.e.} that corresponding to 
$\nabla_{\xi}\ell$) depends on the $2$-jet of the function $f$.(This comes from the fact that 
the Reeb vector field itself depends on the $1$-jet of $f$). More precisely, one has:

\begin{prop}\label{prop:changconnTW}
Let $\theta$ be an adapted contact form and $f$ be a function on $M$. If $\ell$ is a section
of $L$, one has
\begin{equation}\label{eq:changeconnTW} \nabla^{\exp^{2f}\!\theta}\ell - \nabla^{\theta}\ell 
= 2\, \partial_b f 
\otimes\ell + \frac{2}{n+2} \left( df(\xi) + i\Lambda(\nabla^{1,0}\bar\partial_b f) + i(n+1)\,|d_bf|^2
\right)\,\theta\otimes\ell ,
\end{equation}
where $\xi$ is the Reeb field already defined, $\bar\partial_b$ is the $(0,1)$-part of the restriction 
to the contact hyperplane $d_b$ of the differential of a function, $\nabla^{1,0}$ is the $(1,0)$-part
of the Tanaka-Webster connection on $T^*M$ and $\Lambda$ is the 
natural complex trace relative to the K\"ahler-like $2$-form $d\theta$: on a $(1,1)$-form $\beta$,
$$ \Lambda(\beta) = \sum_{\alpha=1}^n \beta(Z_{\alpha},Z_{\bar\alpha}) , $$
where $(Z_{\alpha})$ is an orthonormal basis of $H^{(1,0)}$ with its natural Hermitian metric, and 
$(Z_{\bar\alpha})$ the complex conjugate basis of $H^{(0,1)}$.
\end{prop} 

\smallskip

\begin{proof} Easy computations using the formulas of \cite{dmjc-ricci,cap-slovak,lee-fefferman} 
or \cite[Prop. 2.3]{gover-graham}.
\end{proof}

\medskip

This may look quite strange, as one might expect that the decomposition of $1$-jets of sections of 
$L$ depends on the $1$-jet of the conformal rescaling factor only. However, another choice of connection 
on $L$ leads to a much simpler behaviour: the \emph{Weyl structures}, as defined by  A. \v{C}ap and 
J. Slov\'ak \cite{cap-slovak}, and D. M. J. Calderbank, T. Diemer and V. Sou\v{c}ek \cite{dmjc-ricci}. 

It also turns out that the trivialisation of $\TM$ induced by that choice is exactly the one implicitly
used by previous works on the Cartan bundle such as \cite{cap-survey,gover-graham}: these works indeed
define the Cartan bundle as a bundle which admits the decomposition \eqref{eq:Tfirstdec} in any choice 
of contact form, with explicit transformation properties when changing the contact form. Computations 
given below show that the trivialisations used there are simply those associated with 
Calderbank-Diemer-Sou\v{c}ek's Weyl structures.
   
For simplicity's sake, we shall give here a simple \emph{ad hoc} definition of a Weyl structure, 
that will be enough for our study, and we refer the interested reader to \cite{dmjc-ricci,cap-slovak}
for extra information. Moreover, we
shall concentrate here on  those called \emph{exact} Weyl structures, \emph{i.e.} those associated
with a choice of contact form. Weyl structures are defined in \cite{dmjc-ricci,cap-slovak} in a more 
general way, 
without reference to a specific contact form; the space of Weyl structures is an affine 
space modelled on the space of sections of the cotangent bundle $T^*M$. (This statement is however 
slightly more subtle than it appears, as the basic object is not a covariant derivative on $TM$ 
(nor $K$, nor $L$), but a more algebraic one, and parametrizing the space of Weyl structures with
a given origin is a \emph{non-linear} process, similar to the explicit parametrization of a nilpotent 
Lie group by the exponential map from its Lie algebra; the interested reader is referred to 
\cite{dmjc-ricci,cap-slovak} for more details).

\begin{defi} Let $\eta$ a contact structure on $M$, and $\nabla$ its Tanaka-Webster connection.
The Weyl structure attached to $\eta$ is the linear connection $D$ on $L$ defined by
$D = \nabla + \frac{iR}{2(n+1)(n+2)}\eta$, where $R$ is the Tanaka-Webster scalar curvature
of $\eta$.
\end{defi}

\begin{question}
The definition of Weyl structures in \cite{dmjc-ricci,cap-slovak} is a purely algebraic one in the 
context of parabolic Cartan geometries, and the explicit expression given above for the induced 
connection on $L$ is obtained through the Cartan connection. It would be highly desirable to know 
if there is a way to characterize Weyl connections on $L$ (or, equivalently, on $H$ or $TM$) purely 
in terms of the basic geometric features of CR geometry~? Note however that knowledge on this point 
is not necessary for the construction of the Cartan bundle and its natural connection from basic 
geometric considerations done in the present paper.
\end{question}

\begin{rem}
L. David defined another notion of Weyl structures for CR geometry~\cite{liana-CR}. Both
notions are truly different as L.~David's Weyl structure attached to a choice of contact 
form yields the Tanaka-Webster connection, whereas the notion of Weyl structure used here
does not, as we have seen.
\end{rem}

\smallskip

As before, the choice of a Weyl structure attached to a contact form leads to a global decomposition 
$ TM  = TM/H \oplus H $ induced by the existence of the Reeb vector field, which moreover trivialises 
the $TM/H$-factor. It can be noted that for more general Weyl structures, a generalized Reeb vector 
field exists, \emph{i.e.} a map that sends the trivial line bundle $TM/H$ (but \emph{without} a 
canonical trivialisation) into $TM$, transverse to $H$ \cite{dmjc-ricci,cap-slovak,liana-CR}.
Whenever a contact form $\theta$, its induced Reeb field, and its companion Weyl structure are chosen, 
$$\TM = L \,\oplus\, (H^*)^{(1,0)} \otimes L \,\oplus\, L . $$
If $\ell$ is some holomorphic section of $L$, the corresponding 
element of $\TM$ reads in this decomposition
$$ j^1\ell = \left( \,  \ell \, , \,  D^{1,0}\ell \, ,\, -i \, D_{\xi}\ell \, \right) , $$
or, if we continue to use the Tanaka-Webster connection $\nabla$ attached to $\eta$, 
$$ j^1\ell = \left( \,  \ell \, , \,  \nabla^{1,0}\ell \, , \, - i\,\nabla_{\xi}\ell  + 
\frac{R}{2(n+1)(n+2)}\ell \, \right) .$$
Of course, the factor $(-i)$ in both formulas is by no means necessary. We have introduced it to make
all the computations we will make later fully coherent with the choices made for the Cartan 
vector bundle in \cite{cap-survey,gover-graham}.
 
\begin{notn}
When written with respect to a given contact form, elements of the bundle $\TM$ will be identified 
with triples denoted by $(\ell,\tau,\psi)$.
\end{notn}

The introduction of the correcting factor involving Tanaka-Webster's scalar curvature has the
main effect to cancel out the dependence of the gauge transformation formulae on the second-order
derivative when shifting from a contact form to another.
From Proposition \ref{prop:changconnTW}, a given element $(\ell,\tau,\psi)$ in a given contact
form $\theta$ transforms as follows when the contact form is changed to 
$\hat\theta=\exp^{2f}\theta$:
\begin{equation}\label{eq:changeltT}
\begin{split}
\hat\ell \ = & \ \ell  , \\
\hat\tau \ =  & \ \tau + 2\, \partial_b f\otimes \ell , \\
\hat\psi \ =  & \ \psi - 2\, \iota ({\nabla f})\,\tau - \left( |d_bf|^2 + i\, df(\xi)\right)\ell ,
\end{split}
\end{equation}
where $\iota$ denotes the interior product, and $\partial_b$, $d_b$, and $\xi$ are defined in Proposition
\ref{prop:changconnTW}.
The reader may compare this to \cite[formula (3.2)]{gover-graham} (where $\Upsilon = 2f$) to convince
himself (or herself) that our Cartan bundle is identical to the one obtained by the previously used 
\emph{ad hoc} definition, as the only bundle that is trivialised as 
$$\TM = L \,\oplus\, (H^*)^{(1,0)} \otimes L \,\oplus\, L  $$
in each contact form, with changes of trivialisation given by formulae 
\eqref{eq:changeltT}, see \cite{cap-survey,gover-graham}.

\smallskip

\begin{rem} From this study, it appears that Weyl connections should be more widely used in CR geometry,
as they might provide a more efficient way of treating CR-invariant problems than the usual 
Tanaka-Webster connection. However, since the goal of the present paper is to present an explanation
of the Cartan bundle to an audience that is widely used to computations in the Tanaka-Webster setting, 
we shall continue to do so in the following sections, so that this paper could be compared to previously
existing works.
\end{rem}

\medskip

We now continue our study of the Cartan bundle. This rank $(n+2)$ complex vector bundle has a lot 
of extra structures we shall now detail, the simplest one being a natural determinant:

\begin{prop}
There exists a natural identification of $\Lambda^{n+2}\TM$ with the trivial bundle $M\times\CM$.
\end{prop}
 
\begin{proof}
In a specific Weyl structure, $\TM = L\oplus (H^*)^{(1,0)}\otimes L \oplus (TM/H)^*\otimes L$, so that
$\Lambda^{n+2}\TM = (TM/H)^*\otimes \Lambda^{n}(H^*)^{(1,0)} \otimes L^{n+2} = K\otimes L^{n+2} = \CM$;
Proposition \ref{prop:changconnTW} and the formulas \eqref{eq:changeltT} easily show that this is 
independent of the Weyl structure. \end{proof}

\medskip

\section{Second order geometry and the Fefferman metric}

\smallskip

We shall pursue here our study of the Cartan bundle. The elementary remarks done in the previous
section depend on the first-order geometry of $L$. However, the structure of $\TM$ reveals itself 
only when second-order geometry enters the picture.

\smallskip

Let $\theta,\theta^1,\ldots,\theta^n$ be a local orthonormal coframe such that
$$ d\theta = i \theta^{\alpha}\wedge\theta^{\bar\alpha} $$
(summation convention is understood here, as everywhere in the paper) and let $\ell_{\textrm{ref}}$ be 
a local section of $L$ such that $(\ell_{\textrm{ref}})^{\otimes -(n+2)} = 
\theta\wedge\theta^1\wedge\ldots\wedge\theta^n$. We denote by $\Ric$ the Tanaka-Webster Ricci curvature 
of $\theta$ and by $R$ its scalar curvature (which is one 
half of the usual Riemannian scalar curvature). We will also need the CR Schouten tensor (or the
equivalent real $(1,1)$-form on $H$)
$$\mathcal{P} = \frac{1}{n+2}\left( \Ric - \frac{R}{2(n+1)}\gamma\right),$$ 
(where $\gamma(.,.)=d\theta(.,J.)$ is the metric on $H$); we shall denote with same letter the 
symmetric form and the associated $(1,1)$-form, and things will be 
clear from the context. In the orthonormal coframe above, the coefficients of the $(1,1)$-form are
$\mathcal{P} = iP_{\alpha\bar\beta}\,\theta^{\alpha}\wedge\theta^{\bar\beta}$ (summation convention
intended here again), with 
$$ 
P_{\alpha\bar\beta} = \frac{1}{n+2}\,\left( R_{\alpha\bar\beta} - \frac{R}{2(n+1)} 
\delta_{\alpha\bar\beta}\right) .
$$ 
Another useful object is the Tanaka-Webster torsion, which can be seen as a $J$-antiinvariant 
symmetric bilinear form $\mathcal{A}$ on $H$; hence, it can be written in the local orthonormal 
coframe above as $\mathcal{A}=A_{\alpha\beta}\theta^{\alpha}\otimes\theta^{\beta} + 
A_{\bar\alpha\bar\beta}\theta^{\bar\alpha}\otimes\theta^{\bar\beta}$, with 
$A_{\alpha\beta}$ symmetric in $(\alpha,\beta)$. 

\smallskip

Choose now another local section $\ell$, that we will suppose for a few moments to be 
\emph{CR-holomorphic}, and let $\hat\theta =\exp^{-2\rho}\theta$ be the \emph{volume 
normalized} contact form relative to $\ell$. One has $\exp^{2\rho} = \la \ell ,\ell \ra$ (with respect
to our background section and contact  form
$(\ell_{\textrm{ref}},\theta)$), and for future reference, we note the following useful 
formulas for the function $\rho$ (and subsequent Tanaka-Webster derivatives in the orthonormal coframe 
$\theta^0=\theta,\theta^1,\ldots,\theta^n$). In what follows, a subscript $0,\alpha, \bar\beta$ 
always denotes a derivative with respect to the Tanaka-Webster connection attached to $\theta$. 
If more than one derivative is involved, we adopt the following convention: 
$t_{\gamma\alpha\beta} = \nabla_{Z_{\beta}}\nabla_{Z_{\alpha}}t_{\gamma}$... In case some ambiguity 
might occur, we shall separate the original indices with the extra ones appearing after a derivative by
a comma: $\nabla_{Z_{\beta}}\nabla_{Z_{\alpha}}t_{\gamma}= t_{\gamma,\alpha\beta}$.
First derivatives of $\rho$ are 
\begin{equation}\label{eq:exprhoderivees1}
2\,\rho_{0}\, \exp^{2\rho} = \la \ell_{0}, \ell \ra + \la \ell, \ell_0 \ra \,;\ \
2\,\rho_{\alpha}\, \exp^{2\rho}  = \la  \ell , \ell_{\alpha} \ra \,; \ \ 
2\,\rho_{\bar\alpha}\, \exp^{2\rho}  = \la \ell_{\alpha}, \ell \ra 
\end{equation}
whereas second derivatives are
\begin{equation}\label{eq:exprhoderivees2}
\begin{split}
2\,\rho_{\alpha\beta}\, \exp^{2\rho} \ + \ 4\,\rho_{\alpha}\rho_{\beta}\, \exp^{2\rho} 
& = \la  \ell , \ell_{\alpha\beta} \ra \\
2\,\rho_{\alpha\bar\beta}\, \exp^{2\rho} \ + \ 4\,\rho_{\alpha}\rho_{\bar\beta}\, \exp^{2\rho} 
& =  \la  \ell_{\beta} , \ell_{\alpha} \ra + \la  \ell , \ell_{\alpha\bar\beta} \ra 
= \la  \ell_{\beta} , \ell_{\alpha} \ra - P_{\alpha\bar\beta} \la \ell ,\ell \ra
- \delta_{\alpha\bar\beta} \la \ell ,\psi \ra \\
2\,\rho_{0\alpha}\, \exp^{2\rho} \ + \ 4\,\rho_{\alpha}\rho_{0}\, \exp^{2\rho} 
& = A_{\alpha\beta}\,\la\ell_{\beta},\ell\ra -{\frac{1}{n+2}}\, A_{\alpha\beta,\beta}\,\la\ell,\ell\ra
+ \la\ell_0,\ell_{\alpha}\ra + \la\ell,\ell_{0\alpha}\ra .
\end{split}
\end{equation}
where $\delta_{\alpha\bar\beta}$ is the Kronecker symbol and we recall 
$\psi = -i \ell_0 + \frac{R}{2(n+1)(n+2)}\ell$. 
In all these computations, we have used the fact that $\ell$ is an holomorphic section, including at
second order, so that $2$-jet geometry of $L$ is involved. 

Easy consequences of these formula are the expressions of the curvature and torsion 
of the Tanaka-Webster connection of the volume normalized contact form relative to $\ell$. We shall
first give the explicit expression of torsion $\hat{\mathcal A}$ and scalar 
curvature $\hat R$ (a hat always denotes a quantity attached to $\hat\theta$):
\begin{equation}\label{eq:AetR}
\begin{split}
i\hat A_{\alpha\beta}\, |\ell|^2 \ & = \ \la \ell , \ell_{\alpha\beta} + i A_{\alpha\beta}\,\ell\ra\, ;\\
\ \ \ - \frac{\hat R}{n(n+1)} \ & =  \ \la \ell,\psi\ra + \la\psi ,\ell\ra 
+ \la\ell_{\alpha},\ell_{\alpha}\ra
\end{split}    
\end{equation}
This is easily obtained using formulas (\ref{eq:exprhoderivees1}--\ref{eq:exprhoderivees2}) 
in conjunction with the transformation laws for torsion or the various curvature expressions 
of the Tanaka-Webster connection, as found for instance in \cite[(5.6)-(5.15)]{lee-fefferman}

\begin{rem}
The reader is warned that components of the $\hat\theta$-torsion or curvature elements are given here 
with respect to the $\theta$-orthonormal basis. This convention, which is opposite to that of
\cite{lee-fefferman}, will be in order in the whole paper.
\end{rem}

\smallskip

Note that the expression above for the scalar curvature may give the impression that the Webster
scalar curvature only depends on the $1$-jet of $\ell$. This is actually wrong, as one needs to use
the fact that $\ell_{\bar\beta\alpha}=0$ everywhere to get the expression above for the Webster
scalar curvature. Granted this, it turns out that all second derivatives appearing in the course
of the computations can be related to first derivatives due to the curvature formula:
\begin{equation}\label{eq:curvature}
 - \ell_{\alpha\bar\beta} \ =\ \ell_{\bar\beta\alpha} - \ell_{\alpha\bar\beta} \ =\ 
- i\, \delta_{\alpha\bar\beta}\,\ell_0 \ + \ \frac{1}{n+2}\, R_{\alpha\bar\beta}\,\ell  .
\end{equation}
Formulas \eqref{eq:AetR} depend at most on the second derivative of $\ell$. Their validity can now 
be extended to the jet level. More precisely, the condition for non-zero sections of $L$ to be 
CR-holomorphic is a linear \emph{relation} on their $1$-jets, which defines 
$\TM$ as the sububndle of solutions of this relation. \emph{Prolungating} this relation to $2$-jets 
yields a subbundle $\widetilde{\TM}$ of $J^2L$, and the curvature condition \eqref{eq:curvature} says 
that some part of the $2$-jet is controlled by the $1$-jet (whereas other ``second derivatives'' can be 
freely fixed). On this new subbundle $\widetilde{\TM}$ of $2$-jets, the torsion  
is a function on the open set of $2$-jets which project onto non-zero $0$-jets into $\sym$, the
bundle of symmetric bilinear forms on $H^{(1,0)}$.

\begin{rem}\label{rem:l=0} This function does not extend naturally to the whole bundle $\widetilde{\TM}$.
It is defined only on the subbundle $\widetilde{\TM}^+$ of $2$-jets \emph{whose projection onto $0$-jets
is non-zero},
since there is a factor $|\ell|^2$ in front of the transformation formula for the torsion. But one 
notices that $\hat{\mathcal{A}}\otimes\ell$ does extend as a function into $\sym\otimes L$. Indeed, 
on any section of $\widetilde{\TM}$ that has a non-zero projection onto $0$-jets, 
$$ \hat A_{\alpha\beta}\otimes\ell \ = \ -i\, 
\frac{\la \ell , \ell_{\alpha\beta}\ra \otimes\ell}{|\ell|^{2}} 
+ A_{\alpha\beta}\ = \ -i \,\ell_{\alpha\beta}+ A_{\alpha\beta}\,\ell $$
and this has a well defined meaning for \emph{any} section of $\widetilde{\TM}$.
\end{rem}

\smallskip

The Tanaka-Webster scalar curvature $\hat R$ is another example, but with the notable 
feature that it is defined on the whole bundle $\widetilde{\TM}$ (no tensoring by some
power of $L$ is needed to make it well-defined) 
and it actually factorizes through the natural ``forgetful'' projection 
$j^1:\widetilde{\TM}\subset J^2L \rightarrow \TM\subset J^1L$. Since the whole construction is fully
CR-invariant, we can conclude:

\begin{prop} The Cartan vector bundle $\TM$ is endowed with a natural Hermitian metric of
signature $(1,n+1)$ given by the Tanaka-Webster scalar curvature.
\end{prop}

\smallskip

Of course, in case one considers a section $\sigma$ of $\TM$ which is not the $1$-jet of some 
hoomorphic section of $L$, $\hat R$ is no more equal to the Webster scalar curvature of the volume 
normalized contact form relative to the $0$-jet of $\sigma$. Other interesting examples of the same kind
are the Reeb vector field map 
$$ \hat\xi \ = \ \la\ell,\ell\ra\, \xi - i \la\ell,\ell_{\alpha}\ra\, Z_{\bar\alpha} 
+ i\la\ell_{\bar\alpha},\ell\ra\, Z_{\alpha} $$
(with the summation convention being used here), seen as a function from $\TM$ to $TM$ and transverse
to $H$, or the compatible metric function which assigns to any element of $\TM$ the metric 
$\gamma_{\ell} = d\eta_{\ell}(\cdot, J\cdot)$ on $H$. As before, we do not claim that any of these 
expressions actually represents the Reeb vector field or the compatible metric associated to the
$0$-jet part of a section of $\TM$, as every such expression is to be understood as a function 
on the $1$-jets and not as a $1$-st order differential operator on sections (no differentiation
involved). A last example, which will be of some importance below, is the case of the CR Schouten 
$2$-form. It is easily computed that
$$ \hat P_{\alpha\bar\beta} 
= P_{\alpha\bar\beta} + (\rho_{\alpha\bar\beta} + \rho_{\bar\beta\alpha})
- 2 \rho_{\gamma}\rho_{\bar\gamma}\,\delta_{\alpha\bar\beta} . $$
CR-holomorphicity (at second order) of $\ell$ then implies automatically
$$ \hat P_{\alpha\bar\beta} = -\frac{1}{2}\,\exp^{-2\rho}\,  (\,\la \ell,\psi\ra + \la\psi ,\ell\ra 
+ \la\ell_{\gamma},\ell_{\gamma}\ra \, ) \,\delta_{\alpha\bar\beta}, $$
and using the expression found above for the Tanaka-Webster scalar curvature, we get that the metric
is \emph{pseudo-Einstein} (as is well known, see \cite{lee-einstein}):
$$ \hat{\mathcal{P}} = \frac{\hat R}{2n(n+1)}  \,\hat\gamma . $$
\medskip

As another application of our definition of $\TM$ and our use of Weyl structures attached to a given 
gauge, we will recall the well-known definition of the Fefferman metric on $C=L/\RM_+^*$ attached to 
an element of $\TM$, see \cite{farris,lee-fefferman}. 

We begin first by recalling the construction of the Fefferman on $C$ associated to a choice of a 
compatible contact form $\theta $. This endows $L$ with a Hermitian metric and identifies $C$ to the 
circle bundle of unit elements in $L$. The Weyl structure attached to $\theta$ is a metric
connection on $L$, hence it yields a $S^1$-invariant connection $1$-form $\varpi$ with values in $i\RM$ 
on the tangent bundle to $C$. Letting $\pi : C \rightarrow M$ be the natural bundle projection, the 
Fefferman metric attached to $\eta$ is
\begin{equation}\label{eq:Fef1}
g_{\mathrm{F}} = i\varpi \odot \pi^*\theta + \frac{1}{2} \,\pi^*\gamma 
\end{equation}
where $\odot$ is the symmetrized tensor product for $1$-forms $u\odot v = u\otimes v + v\otimes u$, 
and we have denoted $\gamma = d\theta (\cdot,J\cdot)$.
The Fefferman metric is a natural bridge between CR and conformal geometries: if $\hat\theta = 
\exp^{-2u}\theta$, then $\hat g_F = \exp^{-2u} g_F$. For instance, it is expected
that local CR invariants of $M$ can be reconstructed from conformal invariants on the 
$(2n+2)$-dimensional manifold $C$. This remark will be of importance below.

\begin{rem} The Fefferman metric lives here on the circle bundle associated to $L$, not $L^*$,
nor $K$. This explains the sign discrepancies with \cite{farris,lee-fefferman}.
\end{rem}

\smallskip

To pursue our construction later on, we shall need the expression of the Ricci curvature of the Fefferman
metric attached to a contact form $\theta$. We need first a bit of notation: set
\begin{equation}
%\begin{split}
\mathcal{T} = \frac{1}{n+2}\left(\frac{1}{2(n+1)} d_bR + i\,\div{\mathcal A}\right) ,\ \ 
S = \frac{1}{n}\left( \div\mathcal{T} - |P|^2 + |A|^2 \right) , 
%\end{split}
\end{equation}
where $\div$ is the Riemannian divergence (with a minus sign in front), and squared norms are
relative to the Hermitian metric on $H^{(1,0)}$, not the Riemannian metric on $H$. 
In the orthonormal coframe 
$\theta^0=\theta,\theta^1,\ldots,\theta^n$ satisfying $d\theta=i\theta^{\alpha}\wedge\theta^{\bar\alpha}$
as before, the $(1,0)$-part of $\mathcal{T}$ reads
$$ \ T_{\alpha} = \frac{1}{n+2}\left( \frac{R_{,\alpha}}{2(n+1)} - i A_{\alpha\beta,\bar\beta}\right) $$
and
$$ S = -\frac{1}{n}\left( T_{\alpha,\bar\alpha} + T_{\bar\alpha,\alpha} + P_{\alpha\bar\beta}
P_{\bar\alpha\beta} - A_{\alpha\beta}A_{\bar\alpha\bar\beta}\right) $$
(the summation convention being once again used for pairs of repeated indices).
The Fefferman's metric Ricci curvature has been computed by J. Lee \cite{lee-fefferman}, it reads:
\begin{equation}
\label{eq:ricFefsimple}
\Ric_{\mathrm{F}} = \frac{R}{n+1}\, g_{\mathrm{F}}  - 2n\,\varpi^2 - 2n\, S\,\theta^2
 + n\,(\mathcal{P} + \mathcal{A}^J) + n \,\mathcal{T}^J\odot\theta \ ,
\end{equation}
where $\varsigma^2=\varsigma\otimes\varsigma$ for a $1$-form $\varsigma$, the product $\odot$ 
has already been defined, and
$${\mathcal A}^J(\cdot ,\cdot)={\mathcal A}(J\cdot ,\cdot), \ \, \mathcal{T}^J = \mathcal{T}\circ J.$$
Equivalently, in the local coframe $(\varpi,\theta,\theta^{\alpha},\theta^{\bar\alpha})$ lifted to $C$,
\begin{equation}
\label{eq:ricFef}
\begin{split}
\frac{1}{n}\,\Ric_{\mathrm{F}} \ = \ & - \frac{R}{n(n+1)}\,i\varpi\odot\theta 
- 2\,(\varpi^2 + S\,\theta^2)  \\ 
& + \left( P_{\alpha\bar\beta} + \frac{R}{2n(n+1)}\,\delta_{\alpha\bar\beta}\right)
\,\theta^{\alpha}\odot\theta^{\bar\beta}  \\
&  + i\,\left( A_{\alpha\beta}\,\theta^{\alpha}\otimes\theta^{\beta} - 
A_{\bar\alpha\bar\beta}\,\theta^{\bar\alpha}\otimes\theta^{\bar\beta} + 
T_{\alpha}\,\theta^{\alpha}\odot\theta - T_{\bar\alpha}\,\theta^{\bar\alpha}\odot\theta \right) \ .
\end{split}
\end{equation}
The scalar curvature of the Fefferman metric is easily computed, and the result is
$$ \Scal_{\mathrm{F}} = \frac{2(2n+1)}{n+1} \, R \ .$$ 
Two important remarks are in order here: Fefferman metrics never are Einstein, and 
the Fefferman metric depends on the $2$-jet of the contact form, so that its curvature depends 
on its $4$-jet. It is a remarkable fact however that the transformation formulae for all the elements 
appearing in the Ricci curvature (including $\mathcal{T}$ and $S$) of the Fefferman metric only depend 
on two derivatives of the scaling function (said otherwise, $\mathcal{T}$ and $S$ are
fourth-order only in the direction of a change of the CR structure). In the orthonormal coframe, 
if $\hat\theta = \exp^{-2\rho}\theta$, then
\begin{equation}
\label{eq:transTetS}
\begin{split}
\exp^{-2\rho}\,\hat T_{\alpha} = &\ T_{\alpha}  - i\rho_{0\alpha} - 2 P_{\alpha\bar\beta}\rho_{\beta}
+ 2i \, A_{\alpha\beta}\rho_{\bar\beta}
+ 2 \rho_{\alpha\beta}\rho_{\bar\beta} - 2 \rho_{\alpha\bar\beta}\rho_{\beta}
+ 4 \rho_{\alpha}\rho_{\beta}\rho_{\bar\beta} \ , \\
\exp^{-4\rho}\,\hat S = &\ S - \rho_{00} + 6\, (T_{\alpha}\rho_{\bar\alpha} + 
T_{\bar\alpha}\rho_{\alpha}) 
+ 4i\, ( \rho_{0\bar\alpha}\rho_{\alpha} -  \rho_{0\alpha}\rho_{\bar\alpha}) -(\rho_0)^2 \\
& + 6i\, (A_{\alpha\beta}\,\rho_{\bar\alpha}\rho_{\bar\beta} 
- A_{\bar\alpha\bar\beta}\,\rho_{\alpha}\rho_{\beta})
- 12\, P_{\alpha\bar\beta}\rho_{\bar\alpha}\rho_{\beta} 
+ 4\, (\rho_{\alpha\beta}\rho_{\bar\alpha}\rho_{\bar\beta}+
\rho_{\bar\alpha\bar\beta}\rho_{\alpha}\rho_{\beta})\\
&  + 8\, (\rho_{\alpha\bar\beta}+ \rho_{\bar\alpha\beta})\rho_{\bar\alpha}\rho_{\beta} 
+ 12\, (\rho_{\alpha}\rho_{\bar\alpha})^2 .
\end{split}
\end{equation}

\bigskip

\section{The canonical connection}

\smallskip

Our goal here is the construction of the Cartan connection on $\TM$. This mainly involves 
second-order geometry of $L$ and the above computations of Fefferman metric's Ricci curvature.

\smallskip

Recall first that $\TM^+$ (resp. $\widetilde{\TM}^+$) is the bundle of CR-holomorphic $1$-jets (resp.
$2$-jets) that have a \emph{non-zero projection} into $0$-jets. 
The preimage in $2$-jets of an element of $\TM^+$ is an affine subspace 
of $\widetilde{\TM}^+$. Its volume normalized contact form is well defined, and so is its
Fefferman metric. According to what have been said above, the Ricci curvature of this Fefferman 
metric \emph{does make sense} here (although the whole curvature does not). 

We shall then say that the Ricci curvature of a Fefferman metric is \emph{maximally constant} if it
is the closest possible to that 
of the Fefferman metric of a model space for CR geometry, \emph{i.e} a CR manifold admitting a 
pseudo-Einstein contact form with vanishing torsion. 
Looking back to formulas (\ref{eq:ricFefsimple}-\ref{eq:ricFef}), this means 
that all terms are zero, or constants depending only on $R$ and universal constants in the dimension. 
More precisely, a contact form has a Fefferman metric with maximally constant Ricci curvature iff.
$$ \mathcal{P} = \frac{R}{2n(n+1)}\,\gamma ,\ \ \mathcal{A} = 0, \ \ \mathcal{T} = 0,\ \
\textrm{and } \, S = -\frac{R^2}{4n^2(n+1)^2}\, .$$ 

\smallskip

\begin{ex}[the standard CR sphere in dimension $3$]
The standard left-invariant $1$-forms on $\cerc^3$ are $(\sigma_0,\sigma_1,\sigma_2)$ and satisfy
$d\sigma_i = 2\,\sigma_{i+1}\wedge\sigma_{i+2}$.
Choosing $$ \theta = \frac{1}{2}\,\sigma_0, 
\ \ \theta^1 = \frac{1}{\sqrt{2}}\,(\sigma_1 + i\sigma_2)\, ,$$
yields the Tanaka-Webster connection $1$-form $\omega_1^1 = - 4i\,\theta$, torsion $A_{11} = 0$,
and curvature $R=4$.
The Weyl covariant derivative on $L$ has then a connection $1$-form $\varpi_1^1 = -i\,\theta$
in the frame induced by the choice of $\theta\wedge\theta^1$ as background section on the canonical 
bundle $K$. The Fefferman metric lives on $C=\cerc^1\times\cerc^3$; at a point $(v,q)$, the connection
form is $ \varpi = i dv + \varpi_1^1$, so that
\begin{equation}\begin{split} 
g_{\mathrm{F}} \ & = 
(\theta - dv )\odot\theta + \frac{1}{2}\left( (\sigma_1)^2 + (\sigma_2)^2 \right)\\
& = \frac{1}{2}\left( (\sigma_0)^2 + (\sigma_1)^2 + (\sigma_2)^2 
- \,\sigma_0\odot dv \right) \, .
\end{split}
\end{equation}
This is nothing else but the standard \emph{Einstein universe} metric, \emph{i.e.} the product 
Lorentzian metric on $\cerc^1\times\cerc^3$. It is the model space for conformally 
flat Lorentzian geometry in dimension $4$, but it is here seen in a disguise. Indeed, letting
\begin{equation*}\begin{split} 
\Phi : \cerc^1\times\cerc^3 & \longrightarrow \,\cerc^1\times\cerc^3\\
(v,q) & \longmapsto (v, v^{-1}q) 
\end{split}
\end{equation*}
where the element $v$ of $\cerc^1\subset\CM$ acts on an element $q=(z_1,z_2)$ of $\cerc^3\subset\CM^2$
by $v^{-1}q = (v^{-1}z_1,v^{-1}z_2)$. 
Then one can easily check that
$$ \frac{1}{2}\,\Phi^*(-dt^2 + g_{\cerc^3}) = g_{\mathrm{F}}\, ,$$
thus substantiating our claim.
\end{ex}

\smallskip

Our main result can now be stated as follows.

\begin{theo}\label{th:embed} There is a natural linear map  
\begin{equation*}\label{eq:r} r : \ \TM \longrightarrow J^1\TM \cap J^2L \, .\end{equation*}
The map $r$ embeds $\TM^+$ in $J^2L$ as the bundle of $2$-jets of CR-holomorphic sections of $L$ with 
pointwise Fefferman metric of maximally constant Ricci curvature 
with constant scalar curvature. It is a section of the 
projection $J^1\TM \rightarrow\TM$, so that the image of $r$ is transverse to the principal 
part $T^*M\otimes\TM$ in $J^1\TM$.
\end{theo} 

\begin{rem} It seems unfortunately impossible to give a sense to the \emph{maximally constant Ricci
curvature} equations 
at the jet level for elements outside $\widetilde{\TM}$ (\emph{i.e.} for jets having vanishing 
projection in $0$-jets), see the proof below for the explicit formulas. The image of $r$ is however 
given by linear equations in $J^1(J^1L)$, that are equivalent to the maximally constant Ricci curvature 
equations with constant scalar curvature on $\widetilde{\TM}$, see the proof below.
\end{rem}

\smallskip

Given Theorem \ref{th:embed}, one can now conclude:

\begin{defi} 
The canonical Cartan connection is the natural linear connection provided by the embedding $r(\TM)$ 
in $J^1\TM$ transverse to the principal part.
\end{defi}

\begin{rem}\label{rem:conforme} 
This should be compared to the conformal case, where the Cartan bundle is the sub-bundle
of $2$-jets of the line bundle of $1$-densities 
that are pointwise Einstein, and the Cartan connection is 
then automatically deduced from that very definition of the Cartan bundle \cite{beg,pg-cartan}. 
In the CR 
case, things are at the same time a bit simpler since part of the condition is of first-order (and 
translates as a holomorphicity condition on the sections of $L$), and more complicated since the
equation to be imposed to the Fefferman metric is not the Einstein condition but a translation
thereof. This explains why the Cartan bundle in the CR case can be described at
the level of $1$-jets of sections of $L$, but, as we have already seen, its full geometry is
revealed only when considering $2$-jets. (We have already met this phenomenon when describing the 
natural Hermitian product on $\TM$).

\smallskip

Comparison with the conformal case shows that an alternative definition of the Cartan bundle could have
been: the bundle of CR-holomorphic $2$-jets of $L$ whose Fefferman metric has pointwise maximally 
constant Ricci curvature
with constant scalar curvature. The main theorem of this section shows that both definitions 
coincide ; however, the definition of the Cartan connection in this last approach would have involved 
jets of third order, which makes computations much less tractable. Moreover, the use of third-order
derivatives may involve some non-trivial translations in the jet bundle $J^3L$. An example thereof
is the conformal case where the Einstein condition implies the vanishing of the Cotton tensor, so 
that the sought connection is not in the intersection of the holonomic lifts with the prolungated
lifts (see the proof below) 
but in some translated affine space, due to the differential Bianchi identity \cite{beg,pg-cartan}. 
This phenomenon
is purely third-order and is not to be met below.

\smallskip

We also note another 
surprising feature of the construction. As the Cartan bundle for \emph{conformal} geometry 
is a bundle of pointwise Einstein jets of metrics \cite{beg,pg-cartan},
and the conformal geometry of the Fefferman space
reflects the CR geometry of the basis manifold, one could have expected the \emph{CR Cartan bundle} to be
a bundle of jets of sections of $L$ giving rise to \emph{Einstein} Fefferman metrics. As already 
observed, this is impossible. The geometry of the CR Cartan bundle and connection reflects this fact 
as it forces parallel sections of $\TM$ to induce a Fefferman metric whose Ricci tensor looks like that
of the Einstein universe and not a Lorentzian 
space of constant curvature. This implies that, if one wishes to
interpret the conformal Cartan bundle of the Fefferman metric as the pull-back of the CR Cartan bundle 
of the basis, some translation in needed to make contact between Einstein universe-like and Einstein
jets of metrics. The author intends to come back on this point in a future paper. 
\end{rem}

\medskip

{{\noindent\it Proof of Theorem \ref{th:embed}}}. -- 
This is a long one, which involves many computations. We begin with a few general considerations.
Given any connection on $L$, an element $\sigma$ of $J^1L$ can be written as $\sigma=(\lambda,\ell)$ with 
$\lambda$ in $T^*M\otimes L$ and $\ell$ in $L$. Choosing moreover an auxiliary connection on $TM$, an 
element $\Sigma$ of $J^1(J^1L)$ is then $\Sigma =(\Lambda,\lambda_1;\lambda,\ell)$, with $\Lambda$ in
$T^*M\otimes T^*M\otimes L$, $\lambda_1$ in $T^*M\otimes L$, and $\lambda,\ell$ as above, whereas an
element of $J^2L$ is of the form $(\Phi,\lambda,\ell)$, with $\Phi$ in $Sym^2T^*M\otimes L$.

Denoting by the same letter $D$ the connection on $L$ and the auxiliary connection on $TM$, one has
the general curvature identity on sections of $L$
$$ (D^2\ell)_{U,V} - (D^2\ell)_{V,U} = R^D_{U,V}(\ell) - D_{T^D_{U,V}}\ell \ ,$$ 
where $R^D$ is the curvature of $L$ and $T^D$ the torsion on $TM$. This gives the natural embedding of 
$J^2L$ in $J^1(J^1L)$:
\begin{equation}
\label{eq:nat-emb}
(\Phi,\lambda,\ell) \longmapsto (\Lambda = \Phi + \frac{1}{2}\, R^D(\ell) -\frac{1}{2}\, 
\lambda\circ T^D\, ,
\lambda \, ;\lambda \, ,\ell) 
\end{equation}
If $TM$ is further split into $TM/H\oplus H$, then for any $\theta$ such that $H=\ker\theta$, one has in 
$J^1L$,
$$\sigma = (\theta\otimes\psi + \tau , \ell ) $$
with $\tau$ in $H^*\otimes L$ and $\psi$ in $L$, and in $J^1(J^1 L)$,
$$ \Sigma = (\theta\otimes\theta\otimes\phi + \theta\otimes m + \varphi\otimes\theta + \mu\, , 
\theta\otimes\psi_1 + \tau_1\,; \theta\otimes\psi + \tau \, , \ell ) $$
with $\phi$ and $\psi_1$ in $L$, $\tau_1$, $\varphi$ and $m$ in $H^*\otimes L$, and $\mu$ in 
$H^*\otimes H^*\otimes L$.
Suppose now that the connection on $TM$ preserves $\theta$ and a vector field $\xi$ transverse
to $H$ such that $\theta(\xi)=1$ (thus it preserves the splitting; examples are the Tanaka-Webster
connection and the \emph{non-Ricci-corrected} Weyl connection in the language of 
\cite{dmjc-ricci,cap-slovak}). 
Then the natural embedding \eqref{eq:nat-emb} specializes in the straightforward way: 
\begin{equation}\label{eq:J2dansJ1J1}
\begin{split}
\phi = &\ \Phi(\xi,\xi) ,\\
\varphi(X) = &\ \Phi(X,\xi) + \frac{1}{2}\, R^D_{X,\xi}(\ell) + \frac{1}{2}\,\tau(T^D_{\xi,X}) , \\
m(X) = &\ \Phi(X,\xi) - \frac{1}{2}\, R^D_{X,\xi}(\ell) - \frac{1}{2}\,\tau(T^D_{\xi,X}),\\
\mu(X,Y) = &\ \Phi(X,Y) + \frac{1}{2}\, R^D_{X,Y}(\ell) - \frac{1}{2}\,d\theta(X,Y)\,\psi ,
\end{split}
\end{equation}
(note that we have replaced $T^D_{X,Y}$ by its value in the last line). 

If the connection on $L$ is of type $(0,1)$ (as, \emph{e.g.}, Tanaka-Webster or Weyl connections), then
$\sigma$ as above is in $\TM$ iff. $\tau$ belongs to $(H^*)^{(1,0)}$. If moreover the auxiliary 
connection
preserves the complex structure on $H$ (as do both Tanaka-Webster or Weyl connections), then
$ \Sigma$ as above is in $J^1\TM$ iff.
$$ \tau, m \in (H^*)^{(1,0)}\otimes L \ \textrm{ and }\ \mu \in H^*\otimes (H^*)^{(1,0)}\otimes L . $$

\smallskip

Now, choose an element $\sigma = (\theta\otimes\psi + \tau , \ell)$ in $\TM$. The affine spaces
\begin{equation}
\begin{split}
P_{\sigma} & = \{ \sigma_1 \in J^1\TM \ , \ j^0\sigma_1 = \sigma\} \ \ \textrm{(CR-holomorphic 
prolungated lifts)}\\
H_{\sigma} & = \{ \sigma_2 \in J^2L \ , \ j^1\sigma_2 = \sigma\} \ \ \textrm{(holonomic lifts)} 
\end{split}
\end{equation}
have an intersection in $J^1(J^1L)$ whose elements are of the following type
\begin{equation}\label{eq:quadru}
 (\theta\otimes\theta\otimes\phi + \theta\otimes m + \varphi\otimes\theta  + \mu\, ,
\ \theta\otimes\psi + \tau\, ;\ 
\theta\otimes\psi + \tau\, ,\ \ell ) 
\end{equation}
and satisfy moreover
\begin{equation}\label{eq:definissant}\begin{split}
& \tau, m \in (H^*)^{(1,0)}\otimes L , \ \ \mu \in H^*\otimes (H^*)^{(1,0)}\otimes L ,\\
& \varphi(\bar Z) = - R^D_{\xi,\bar Z}\ell + \tau\circ T^D_{\xi,\bar Z}, \ 
\textrm{ for any } Z \in H^{(1,0)},\\
& \varphi(Z) = m(Z) - R^D_{\xi,Z}(\ell) + \tau\circ T^D_{\xi,Z}, 
\textrm{ for any } Z \in H^{(1,0)},\\
& \mu(\bar Z, U) = R^D_{\bar Z,U}(\ell) - \frac{1}{2}\,d\theta(\bar Z,U)\psi,\ 
\textrm{ for any } Z \in H^{(1,0)} \textrm{ and any } U \in H.
\end{split}
\end{equation}
This shows that lots of components of the elements of $P_{\sigma}\cap H_{\sigma}$ are already 
determined by this intersection. If we want to find a unique element of $P_{\sigma}\cap H_{\sigma}$,
the only things that are not yet fixed are the components $\varphi^{(1,0)}$ (or $m^{(1,0)}$), 
$\mu^{(1,0),(1,0)}$ and $\mu(\xi,\xi)$.

\smallskip

Fixing the remaining components will now involve the ``maximally constant Ricci curvature with constant 
scalar curvature'' condition. For the sequel of the proof, we now specialize the connection to be the 
usual  Tanaka-Webster derivative. This enables us to work with classical transformation formulas of CR 
geometry. For notational simplicity, we shall also keep notations such as $\ell_{\alpha}$, 
$\ell_{\alpha\bar\beta}$, etc... for these Tanaka-Webster derivatives, rather than 
abstract notations as $m$, $\varphi$ or $\mu$ (note that it has the advantage of clearly making a 
difference between the actual choice of the Tanaka-Webster connection and the previous general
considerations). Straightforward calculations show that conditions (\ref{eq:quadru}-\ref{eq:definissant}), 
\emph{i.e.} prolungation of the CR-holomorphicity condition seen as the intersection 
$P_{\sigma}\cap H_{\sigma}$, are equivalent to:
\begin{equation}
\label{eq:outils}
\begin{cases}
\ell_{\bar\alpha} = 0,\ \ell_{\bar\alpha\beta}=0,\ \ell_{\bar\alpha 0}=0  
\, ,\\ 
\ell_{\alpha\bar\beta} = - P_{\alpha\bar\beta}\,\ell - \delta_{\alpha\bar\beta}\,\psi \, ,\\
\ell_{0\bar\alpha} = - \frac{1}{n+2}\, A_{\bar\alpha\bar\beta,\beta}\,\ell 
+ A_{\bar\alpha\bar\beta}\,\ell_{\beta} 
\end{cases}
\end{equation}
(note that we shall not be completely coherent here as we shall mix a little bit our Tanaka-Webster 
choice with some Weyl derivative, by using $\psi = -i \ell_0 + \frac{R}{2(n+1)(n+2)}\,\ell$ rather than 
$\ell_0$; we hope that the remarks done in the previous sections convinced the reader that it is 
a more natural choice).
These conditions precisely define a subset of $J^1(J^1L)$ that is canonically
identified to $\widetilde{\TM}$ in $J^2L$.
\smallskip

We will now determine the remaining components $\ell_{\alpha\beta}$, $\ell_{0\alpha}$, and $\ell_{00}$
for an elements of $J^2L$ having non-zero $0$-jet. As already observed, the Fefferman metric determined 
by such an element has maximally constant Ricci curvature 
iff. $A_{\alpha\beta}$ and $T_{\alpha}$ vanish and $S$ is completely determined by $R$. 
Whenever $\ell$ is non-zero, \eqref{eq:AetR} shows that the Tanaka-Webster torsion of the
volume normalized ($2$-jet of) contact form $\hat\theta = \exp^{-2\rho}\theta$ being zero 
is equivalent to
\begin{equation}\label{eq:tordefi}
  \ell_{\alpha\beta} + i A_{\alpha\beta}\,\ell = 0 . 
\end{equation}
We can extend this everywhere, so that $\ell_{\alpha\beta}$ is fully determined:
from Remark \ref{rem:l=0} on torsion as a function on the space of $2$-jets in section 4, one
sees that this is the same as requiring that the linear function $\hat{\mathcal{A}}\otimes\ell$
vanishes.

\smallskip

We now manage the condition on $T_{\alpha}$. For this we replace in \eqref{eq:transTetS} every
occurrence of a derivative of $\rho$ by its expression in terms of the $2$-jet $\sigma$, taken
from equations (\ref{eq:exprhoderivees1}-\ref{eq:exprhoderivees2}). 
We will also make free use of the conditions we have 
found on the $2$-jet that comes from the intersection of $P_{\sigma}$ and $H_{\sigma}$, \emph{i.e.} 
formulas \eqref{eq:outils}. Said otherwise, we now consider only sections of $\widetilde{\TM}$ rather
than $J^1(J^1L)$ and funtions defined on this reduced bundle, so that relations \eqref{eq:outils}
are definitely ``plugged in'' the construction. On the contrary, we shall keep aside the 
torsion-free condition \eqref{eq:tordefi}.

\smallskip

Our goal is to show that $\ell_{0\alpha}$ is determined by the vanishing of $\hat T_{\alpha}$
whenever $\ell$ is non-zero. More precisely, we shall prove that, given \eqref{eq:outils}, this is 
equivalent to \eqref{eq:tordefi} together with 
\begin{equation}
\label{eq:l_0alpha} 
 -i\,\ell_{0\alpha} + \frac{i}{2(n+1)(n+2)}\, R_{,\alpha}\,\ell + \frac{i}{2(n+1)(n+2)}
\, R\,\ell_{\alpha} - P_{\alpha\bar\beta}\,\ell_{\beta} + T_{\alpha}\,\ell = 0
\end{equation}
(remember that $R_{,\alpha}$ denotes the derivatives of $R$ in the coframe above). 
Using $\psi = -i \ell_0 + \frac{R}{2(n+1)(n+2)}\,\ell$ as in section 3, \eqref{eq:l_0alpha} 
is equivalent to:
\begin{equation}
\label{eq:psi_alpha} 
\psi_{\alpha} = P_{\alpha\bar\beta}\,\ell_{\beta} - T_{\alpha}\,\ell \ . 
\end{equation}
To do so, one writes the transformation law \eqref{eq:transTetS} for $T_{\alpha}$ in the following form 
\begin{equation*}
\begin{split}
\hat T_{\alpha} 
= &\ \left( \exp^{2\rho}\, T_{\alpha}  - i\,\rho_{0\alpha}\exp^{2\rho} 
- P_{\alpha\bar\beta}\,\rho_{\beta}\exp^{2\rho}
+ i \, A_{\alpha\beta}\,\rho_{\bar\beta}\,\exp^{2\rho} \right) 
- \left( 2\,\rho_{\alpha\bar\beta}\,\exp^{2\rho} + P_{\alpha\bar\beta}\,\exp^{2\rho}\right)\rho_{\beta}\\
&\ + \left( i \, A_{\alpha\beta}\,\exp^{2\rho} + 2\,\rho_{\alpha\beta}\,\exp^{2\rho}
+ 4\,\rho_{\alpha}\rho_{\beta}\,\exp^{2\rho}\right)\rho_{\bar\beta} ,
\end{split}
\end{equation*}
In the previous formula, the second set of parentheses satisfies
$$ - \left(2\,\rho_{\alpha\bar\beta}\,\exp^{2\rho} + P_{\alpha\bar\beta}\,\exp^{2\rho} \right) 
\ = \ \la\ell ,\psi\ra \,\delta_{\alpha\bar\beta}
$$
as a result of prolungation of the CR-holomorphicity condition, and more precisely of the third formula
in \eqref{eq:outils} (the metric is pseudo-Einstein, see section 4 above), so that this term is always 
zero on $\widetilde{\TM}$. For the last set of parentheses, we get
$$
\left( i \, A_{\alpha\beta}\,\exp^{2\rho} + 2\,\rho_{\alpha\beta}\,\exp^{2\rho}
+ 4\,\rho_{\alpha}\rho_{\beta}\,\exp^{2\rho}\right)\rho_{\bar\beta} = 
\la \ell, A_{\alpha\beta}\,\ell - i \ell_{\alpha\beta} \ra \, \rho_{\bar\beta} ,
$$
hence, this yields
\begin{equation*}
\hat T_{\alpha} = \exp^{2\rho} 
\left( T_{\alpha}  - i\rho_{0\alpha} - P_{\alpha\bar\beta}\rho_{\beta}
+ i \, A_{\alpha\beta}\rho_{\bar\beta} + \exp^{-2\rho}\la\ell ,\psi\ra \,\rho_{\alpha} \right)
+  \la \ell, A_{\alpha\beta}\,\ell - i \ell_{\alpha\beta} \ra \, \rho_{\bar\beta} .
\end{equation*}
We now replace $\rho_{0\alpha}$ by its expression in terms of first and second derivatives of 
$\ell$. Letting $w=\psi_{\alpha} - P_{\alpha\bar\beta}\,\ell_{\beta} - T_{\alpha}\,\ell$, a few lines of
computation, using the expression for $\rho_{0\alpha}$ found above, give:
\begin{equation*}
\begin{split}
\hat T_{\alpha} \ = \ & \frac{i}{2}\left( 4\exp^{2\rho}\rho_0\rho_{\alpha}
- \la\ell_0,\ell_{\alpha}\ra - \exp^{-2\rho}\la\ell,\ell_{\alpha}\ra\la\ell,\ell_{0}\ra\right)\\
& \, - \, \frac{i}{2}\,\la \ell, w\ra \, + \, \la \ell, A_{\alpha\beta}\,\ell - i \ell_{\alpha\beta} \ra 
\,\rho_{\bar\beta}
\end{split}
\end{equation*}
It is always true that $4\exp^{4\rho}\rho_0\rho_{\alpha} = \la\ell,\ell_{\alpha}\ra\,
\left( \la\ell_0,\ell\ra + \la\ell,\ell_0\ra\right)$, and, moreover, one notices that 
$\la\ell,\ell_{\alpha}\ra\,\la\ell_0,\ell\ra = \la \ell_0,\ell_{\alpha}\ra$ as soon as
$\ell$ is non-zero, so that 
\begin{equation*}
\hat T_{\alpha} \ =  
\ - \, \frac{i}{2}\,\la \ell, w\ra \, + \, \la \ell, A_{\alpha\beta}\,\ell - i \ell_{\alpha\beta} \ra 
\,\rho_{\bar\beta} .
\end{equation*}
One can then conclude that the equations $\hat{\mathcal T}=0$ and $\hat{\mathcal A}\otimes\ell=0$ (at 
the level of jets) determine $\psi_{\alpha}$ whenever $\ell$ is non-zero.

Note here that one cannot extend this outside the set of (jet-)solutions of the torsion equation
\eqref{eq:tordefi}; indeed, whenever $\ell$ is non-zero one has:
\begin{equation*}
\hat T_{\alpha} \ =  \ - \frac{i}{2}\,\la \ell, w\ra \, 
+ \frac{1}{2}\, \la \ell, A_{\alpha\beta}\,\ell 
- i \ell_{\alpha\beta} \ra \,\frac{\la\ell_{\beta},\ell\ra}{|\ell|^2} \ 
= \ - \frac{i}{2}\,\la \ell, w\ra \, 
+ \frac{1}{2}\, \la \ell_{\beta}, A_{\alpha\beta}\,\ell - i \ell_{\alpha\beta} \ra 
\end{equation*}
or, equivalently,
\begin{equation*}
\hat T_{\alpha}\otimes\ell \ = - \frac{i}{2}\, w \, + 
\,\frac{\la\ell_{\beta},A_{\alpha\beta}\ell - i \ell_{\alpha\beta}\ra}{|\ell|^2}
\end{equation*}
and the right hand side makes no sense on the zero-set of $\ell$; however the function 
$\hat T_{\alpha}\otimes\ell$ has a well-defined meaning on the set of solutions of
\eqref{eq:tordefi} in $H_{\sigma}\cap P_{\sigma}$.%\mh{check again !}

\smallskip

It remains to consider the case of $\hat S$ and the condition on scalar curvature. For the latter
one, it is well-known that the scalar curvature of the Fefferman metric of
$\hat\theta = \exp^{-2\rho}\theta$ equals a constant multiple
of the Tanaka-Webster scalar curvature $\hat R$. Easy computations show that relations found so far imply
that the scalar curvature of the Fefferman metric attached to an element of $\widetilde{\TM}^+$ has
vanishing derivatives in the direction of the contact hyperplane $H$. Moreover, the derivative in
the direction of the Reeb field is given by
$$ - \frac{\hat{R}_{,0}}{n(n+1)} \ = \ \la \ell , y \ra + \la y , \ell \ra \ = 2\,\Re (\la \ell , y\ra ) 
\ $$
where we have denoted
$$ y = \psi_0 + i S\,\ell + 2i\, T_{\bar\alpha}\,\ell_{\alpha} + \frac{i\, R}{2(n+1)(n+2)}\,\psi .$$
The computation for $\hat S$ is much longer but similar to the previous ones and is left to the reader. 
Taking into account all relations found so far, one gets at the end:
$$ |\ell|^{-2}\left(\hat{S} + \frac{4\, {\hat R}^2}{n^2(n+1)^2}\right) \ = \ -\frac{i}{2}
\left( \la \ell , y \ra  - \la y , \ell \ra \right) = \Im (\la \ell , y\ra ) . $$
As a result, $\la \ell, y\ra = 0$ and $y=0$ whenever $\ell$ is non-zero. This determines the last
free component in the $2$-jet of $\ell$ and concludes the proof.
\endproof

\smallskip

We can now resume the definition of the natural Cartan connection $\cnabla$ in a trivialisation of the
Cartan bundle $\TM$ induced by a choice of local basis for $L$, its companion volume-normamlised 
contact form $\theta$, and its Weyl structure. If $(\ell,\tau,\psi)$ is a section of $\TM$ and 
$\nabla$ is the Tanaka-Webster connection of $\theta$, the Cartan derivative  
$\cnabla$ is defined as follows. For a vector $Z$ of type $(1,0)$ in $H$,
\begin{equation}
\label{eq:Cartanconn1}
\cnabla_Z\begin{pmatrix}\ell \\ \tau\\ \psi\end{pmatrix} \ =
\ \begin{pmatrix} 
 \nabla_Z\ell  - \tau(Z)\\ 
 \nabla_Z\tau + i\,{\mathcal A}(Z,\cdot)\otimes\ell \\ 
 \nabla_Z\psi - \tau\circ{\mathcal P}(Z) + {\mathcal T}(Z)\,\ell 
\end{pmatrix}\ ,
\end{equation} 
whereas for a vector $\bar Z$ of type $(0,1)$ in $H$,
\begin{equation}
\label{eq:Cartanconn2}
\cnabla_{\bar Z}\begin{pmatrix} \ell \\ \tau\\ \psi\end{pmatrix} \ =
\ \begin{pmatrix} 
 \nabla_{\bar Z}\ell   \\
 \nabla_{\bar Z}\tau  + {\mathcal P}(\bar Z,\cdot)\otimes\ell + \gamma(\bar Z,\cdot)\otimes\psi \\
 \nabla_{\bar Z}\psi + i\,\tau\circ{\mathcal A}(\bar Z,\cdot)  + {\mathcal T}(\bar Z)\,\ell
\end{pmatrix}\ ,
\end{equation}
and if $\xi$ is the Reeb field of $\theta$,
\begin{equation}
\label{eq:Cartanconn3}
\cnabla_{\xi}\begin{pmatrix} \ell \\ \tau\\ \psi\end{pmatrix} \ = 
\ \begin{pmatrix} 
 \nabla_{\xi}\ell + \frac{i R}{2(n+1)(n+2)}\,\ell - i\,\psi  \\
 \nabla_{\xi}\tau - i\,\tau\circ{\mathcal P} + \frac{iR}{2(n+1)(n+2)}\,\tau 
    + 2i\,{\mathcal T}^{(1,0)}\otimes\ell \\
 \nabla_{\xi}\psi + \frac{i R}{2(n+1)(n+2)}\,\psi + 2i\,\tau\circ{\mathcal T} + i\,S\,\ell 
\end{pmatrix}\ .
\end{equation}
The interested reader can check that these formulas are identical to the ones given in previous
works such as \cite{gover-graham}.

\medskip

\begin{cor}
The Cartan connection is a Hermitian connection which respects the natural determinant of $\TM$.
\end{cor}

\begin{proof} Easy computations; remember that the Hermitian product on $\TM$ is given by scalar 
curvature, so that the constant scalar curvature condition in the definition of the Cartan connection
mandatorily implies the Hermitian character of the connection.
\end{proof}

\medskip

Said otherwise, one may endow the principal $SU(1,n+1)$-bundle of frames of $\TM$ with an
equivariant connection form. Moreover, $\TM$ has a natural isotropic line sub-bundle $\Lambda$, 
whose elements are given by $(\psi,0,0)$
in \emph{any} choice trivialisation of $\TM$ relative to a choice of contact form (in classical 
conformal geometry, these are the \emph{elements at infinity}, see \cite{beg,pg-cartan}). The Cartan 
connection does \emph{not} preserve the line bundle $\Lambda$, as it is easily seen from the
formulas (\ref{eq:Cartanconn1}--\ref{eq:Cartanconn3}). However, the structure group of the frame 
bundle of $\TM$ (without any connection) can be reduced to the parabolic subgroup $P$ of $SU(1,n+1)$ 
that leaves an isotropic line globally invariant. The resulting principal bundle $\mathcal{G}$ is 
the central object in the general theory of parabolic geometries and the connection $\cnabla$ translates
as the canonical Cartan connection on this bundle. 

\bigskip

\section{The curvature of the Cartan connection}

\smallskip

This short section contains no really new results but is included for sake of completeness. 

\smallskip

The curvature of the Cartan connection can be computed from formulas 
(\ref{eq:Cartanconn1}--\ref{eq:Cartanconn3}) with respect
to the Tanaka-Webster connection of some contact form. We first need a bit of notation:
for a given contact form $\theta$ and its associated Tanaka-Webster 
derivative $\nabla$ on $TM$, we let $$ W = R^{\nabla} - \mathcal{P}\kn\gamma .$$ 
It is a $(1,1)$-form on $H$ with values in $(1,1)$-forms, where $\kn$ is the Kulkarni-Nomizu 
product on symmetric bilinear forms; in a compatible orthonormal coframe $(\theta^0=\theta,
\theta^{\alpha},\theta^{\bar\alpha})$ as above,
$$ W_{\alpha\bar\beta\rho\bar\sigma} =  R^{\nabla}_{\alpha\bar\beta\rho\bar\sigma} 
- P_{\alpha\bar\beta}\delta_{\rho\bar\sigma} 
- P_{\rho\bar\sigma}\delta_{\alpha\bar\beta}  
- P_{\alpha\bar\sigma}\delta_{\rho\bar\beta}
- P_{\rho\bar\beta}\delta_{\alpha\bar\sigma}  . $$
It is the Tanaka-Webster equivalent of the Weyl tensor in Riemannian geometry; we shall term
it the Chern-Moser tensor, since it seems to have been first discovered by those two authors 
in \cite{chern-moser}. The second
tensor we will need is third order (here as in all formulas that follow, we shall use a comma to 
separate tensor indices from subsequent derivatives):
$$ V_{\alpha\bar\beta\rho} = A_{\alpha\rho,\bar\beta} + i\,P_{\alpha\bar\beta,\rho}
-i\, T_{\rho}\delta_{\alpha\bar\beta} - 2i\, T_{\alpha}\delta_{\rho\bar\beta}$$
and we shall also need fourth-order tensors 
\begin{equation*}\begin{split}
Q_{\alpha\beta} & = i\, A_{\alpha\beta,0} - 2i\, T_{\alpha,\beta} + 2\, P_{\alpha\bar\rho}
\, A_{\rho\beta} \ , \\
U_{\alpha\bar\beta} & = T_{\alpha,\bar\beta} + T_{\bar\beta,\alpha} 
+ P_{\alpha\bar\rho}P_{\rho\bar\beta} - A_{\alpha\rho}A_{\bar\rho\bar\beta} + 
S\,\delta_{\alpha\bar\beta} \ .
\end{split}
\end{equation*}
The tensor $Q$ is the Cartan tensor discovered in \cite{cartan-dim3,cartan-dim3bis}; 
note that $Q$ is symmetric 
by construction, while $U$ is trace-free by the very definition of $S$. The last piece of 
information is fifth-order:
$$ Y_{\alpha} = T_{\alpha,0} - i\, S_{,\alpha} + 2\,i\, P_{\alpha\bar\rho}\, T_{\rho}
+ 3\, A_{\alpha\rho}T_{\bar\rho} .$$
These are the different components entering the expression of the curvature of the Cartan
connection. Among them, the most important ones are $W$ and $Q$.

\smallskip

The whole curvature now reads as follows on an element $(\ell,\tau,\psi)$ of $\TM$ (summation
convention intended throughout the formulas, as before):
\begin{equation}
\label{eq:courb-Cartan1}
\begin{split}
R^{\cnabla}\begin{pmatrix}\ell\\ \tau\\ \psi\end{pmatrix} \ 
= \ \ & \theta^{\rho}\wedge\theta^{\bar\sigma}\otimes\begin{pmatrix}
0\\
\theta^{\alpha}\otimes\left( i\, V_{\rho\bar\sigma\alpha}\,\ell 
+ W_{\rho\bar\sigma\alpha\bar\beta}\,\tau_{\beta}\right)\\
U_{\rho\bar\sigma}\,\ell
- i\, V_{\bar\sigma\rho\bar\beta}\,\tau_{\beta}
\end{pmatrix}\\
& \  + \ \theta^{\rho}\wedge\theta\otimes 
\begin{pmatrix}
0\\  
\theta^{\alpha}\otimes\left( Q_{\rho\alpha}\,\ell 
+ V_{\rho\bar\beta\alpha}\,\tau_{\beta}\right) \\
Y_{\rho}\,\ell 
-i\, U_{\rho\bar\beta}\,\tau_{\beta}
\end{pmatrix}\\
& \ - \ \theta^{\bar\rho}\wedge\theta \otimes 
\begin{pmatrix} 
0\\
\theta^{\alpha}\otimes\left( i\,U_{\alpha\bar\rho}\,\ell
+ V_{\bar\beta\alpha\bar\rho}\,\tau_{\beta}\right) \\ 
Y_{\bar\rho}\,\ell 
-i\, Q_{\bar\rho\bar\beta}\,\tau_{\beta}
\end{pmatrix}
\end{split}
\end{equation}

\bigskip

We shall now conclude with a few well-known remarks, which we shall give without detailed proofs. 

\smallskip

First of all, a lot of simplifications occur in the especially interesting dimension~$3$. Indeed,
\begin{itemize}
\item $W_{1\bar 1 1\bar 1}=0$ for algebraic reasons;
\smallskip
\item $V_{1\bar 1 1}=0$ and $U_{1 1}=0$ by the very definition of $\mathcal{T}$ and $S$;
\smallskip
\item $Y_1$ is related to the Cartan tensor $Q_{11}$ by the differential Bianchi identity 
for the Cartan connection (this is easily checked by using the formulas for the connection and the
curvature given above).
\end{itemize}
Thus one obtains the well known result that $Q_{11} = 0$ implies that the whole curvature
of the Cartan connection vanishes.

\smallskip

In arbitrary dimension $2n+1$ with $n>1$, things are a bit more complicated since all curvature 
components coexist. However, the Bianchi identity for the Cartan connection implies relations
between them. These relations yield in turn that the vanishing of all
$W_{\rho\bar\sigma\alpha\bar\beta}$'s implies that all of the Cartan curvature components are 
zero. 
In any case one ends with a vanishing curvature for the Cartan connection, and one concludes:

\begin{theo}
Let $(X,H,J)$ be a CR strictly pseudo-convex manifold of dimension $2n+1$, and $\theta$ be 
any choice of compatible
contact form. If $ n>1$, suppose the Chern-Moser tensor $W$ is zero, and if $n=1$, suppose that
the Cartan tensor $Q$ is zero. Then $X$ is \emph{spherical} or \emph{CR-flat}, \emph{i.e.} 
for any point in $x$, there exists a neighbourhood $U_x$ and a compatible contact form $\theta_x$ 
such that $(U_x,H_{|U_x},J_{|U_x},\theta_x)$ is isometric to a neighbourhood of the origin in 
Heisenberg space.
\end{theo}

\begin{proof}
Under the conditions in the Theorem, the whole curvature of the Cartan connection vanishes. 
Hence there exists a wealth of local parallel sections, and one may take an isotropic local 
parallel section $\sigma_x = (\ell,\tau,\psi)$ on a neighbourhood $U_x$ of any point $x$. As it 
is isotropic, its $0$-jet $\ell$ does never vanish so that the volume normalized contact
form $\theta_x$ exists. Moreover, $\sigma_x$ is parallel hence holonomic by construction
of the Cartan connection. As a result, the torsion and CR-Schouten curvature of the associated 
Tanaka-Webster connection vanish, hence its whole curvature, and $\ell^{\otimes [-(n+2)]}$ is a 
closed local section of the CR canonical bundle $K$.

Reasoning as usual, this shows that there exists a local orthonormal coframe 
$\{\theta^{\alpha}_x\}$ on $U_x$ such that
$$ d\theta_x = i\,\theta^{\alpha}_x\wedge\theta^{\bar\alpha}_x \ \textrm{ and } \ 
d\theta^{\beta}_x = 0 \ \ \forall \beta .$$
It remains to show that $U_x$ is isomorphic to to a neighbourhood of the origin in
Heisenberg space. A natural local contactmorphism from $U_x$ to $\RM^{2n+1}$ is then constructed 
as follows: by Poincar\'e lemma, one may let $\theta^{\alpha}
= dx^{\alpha}+i\,dy^{\alpha}$ for some functions $(x^{\alpha},y^{\alpha})$ vanishing at $x$. 
Then 
$$ d\left( \theta_x - \frac12 \sum_{\alpha} x^{\alpha}\,dy^{\alpha} - y^{\alpha}\, dx^{\alpha}
\right) = 0 \, .$$
It then exists a last function $z$ vanishing at $x$ such that
$\theta = dz + \frac12 \sum_{\alpha} x^{\alpha}dy^{\alpha} - y^{\alpha}dx^{\alpha}$.
It is easily seen that (restricting if necessary) these provide coordinates on $U_x$, with
the sought identifications for the contact form and complex structure. This ends the proof.
\end{proof}

\bigskip
 
\begin{rem}
Given a (local) trivialisation of the Cartan bundle, the bundle at infinity $\Lambda$ defines a
natural \emph{development map} from (an open set in) our manifold $M$ into the set of isotropic
lines in $\CM^{n+2}$. In case $M$ is spherical, the Cartan connection is flat, and the Cartan 
bundle may be trivialised (at least locally) by parallel sections. It is (not so easily) shown that the 
(local) embedding above is a CR-equivalence from (an open set in) $M$ and (an open set in) the 
standard CR sphere iff. a frame of parallel sections is used to trivialise $\TM$, in complete 
analogy with what is known for the conformally flat case \cite{pg-cartan}.
\end{rem}

\begin{small}
{\flushleft\sl Acknowledgements}. The author thanks Paul Gauduchon for numerous 
conversations on the subject matter of this work. He is also grateful to Florin Belgun, 
David~M.~J. Calderbank, C.~Robin Graham, and Jan Slov\'ak for their comments on an early version of this paper.
\end{small}

\bigskip

\bibliographystyle{smfplain}

\medskip

\end{document}